\documentclass[12pt,leqno]{article}
\usepackage{amsfonts}
\usepackage{amsmath, amsthm, amsfonts, amssymb, color}
\usepackage{mathrsfs}
\usepackage{color}
\usepackage[notcite,notref]{showkeys}

\theoremstyle{definition}

\newcommand{\scr}[1]{\mathscr #1}
\definecolor{wco}{rgb}{0.5,0.2,0.3}

\numberwithin{equation}{section} \theoremstyle{remark}

\newcommand{\ua}{\uparrow}

\title{{\bf Order-Preservation for Multidimensional Stochastic Functional Differential Equations with Jumps}\footnote{Supported in
 part by  Lab. Math. Com. Sys. and NNSFC(11131003).} }
\author{
{\bf Xing Huang $^{a)}$ and   Feng-Yu Wang$^{a), b)}$}\\
\footnotesize{$^{a)}$School of Mathematical Sciences,
Beijing Normal
University, Beijing 100875, China}\\
 \footnotesize{$^{b)}$Department of Mathematics,
Swansea University, Singleton Park, SA2 8PP, United Kingdom}\\ \footnotesize{wangfy@bnu.edu.cn, F.-Y.Wang@swansea.ac.uk}}
\begin{document}
\allowdisplaybreaks
\def\R{\mathbb R}  \def\ff{\frac} \def\ss{\sqrt} \def\B{\mathbf
B}
\def\N{\mathbb N} \def\kk{\kappa} \def\m{{\bf m}}
\def\dd{\delta} \def\DD{\Delta} \def\vv{\varepsilon} \def\rr{\rho}
\def\<{\langle} \def\>{\rangle} \def\GG{\Gamma} \def\gg{\gamma}
  \def\nn{\nabla} \def\pp{\partial} \def\E{\mathbb E}
\def\d{\text{\rm{d}}} \def\bb{\beta} \def\aa{\alpha} \def\D{\scr D}
  \def\si{\sigma} \def\ess{\text{\rm{ess}}}
\def\beg{\begin} \def\beq{\begin{equation}}  \def\F{\scr F}
\def\Ric{\text{\rm{Ric}}} \def\Hess{\text{\rm{Hess}}}
\def\e{\text{\rm{e}}} \def\ua{\underline a} \def\OO{\Omega}  \def\oo{\omega}
 \def\tt{\tilde} \def\Ric{\text{\rm{Ric}}}
\def\cut{\text{\rm{cut}}} \def\P{\mathbb P} \def\ifn{I_n(f^{\bigotimes n})}
\def\C{\scr C}      \def\aaa{\mathbf{r}}     \def\r{r}
\def\gap{\text{\rm{gap}}} \def\prr{\pi_{{\bf m},\varrho}}  \def\r{\mathbf r}
\def\Z{\mathbb Z} \def\vrr{\varrho} \def\ll{\lambda}
\def\L{\scr L}\def\Tt{\tt} \def\TT{\tt}\def\II{\mathbb I}
\def\i{{\rm in}}\def\Sect{{\rm Sect}}  \def\H{\mathbb H}
\def\M{\scr M}\def\Q{\mathbb Q} \def\texto{\text{o}} \def\LL{\Lambda}
\def\Rank{{\rm Rank}} \def\B{\scr B} \def\i{{\rm i}} \def\HR{\hat{\R}^d}
\def\to{\rightarrow}\def\l{\ell}
\def\8{\infty}\def\I{1}\def\U{\scr U}

\maketitle

\begin{abstract}  Sufficient and necessary conditions are presented for the order-preservation of stochastic functional differential equations on $\R^d$ with non-Lipschitzian coefficients driven by the Brownian motion and Poisson processes. The sufficiency of the conditions extends and improves some  known comparison theorems derived recently for one-dimensional equations and multidimensional equations without delay, and the necessity is new even in these special  situations.   \end{abstract} \noindent
 AMS subject Classification:\  60J75, 47G20, 60G52.   \\
\noindent
 Keywords: Order-preservation, comparison theorem,  stochastic functional differential equation.
 \vskip 2cm

\section{Introduction}

The order-preservation  of stochastic processes is crucial since it enables one to control   complicated processes by using simpler ones. For a large class of diffusion-jump type Markov processes on $\R^d$, the order-preservation property has been well described in the distribution sense (see \cite{CW,W} and references therein), see also \cite{W1} for a study of super processes. To derive the pathwise order-preservation, one establishes the comparison theorem for stochastic differential equations (SDEs) which goes back to \cite{S,Y}.  The study of comparison theorem for one-dimensional SDEs is now very complete, see e.g. \cite{BY, GD, IW, M, O, PY, PZ, YMY} and references therein. Equations considered in these references include forward or backward SDEs with jumps and with delay. The aim of this note is to provide a sharp criterion on the comparison theorem for multidimensional stochastic functional  differential equations (SFDEs) which is   yet unknown  in the literature.

Throughout the paper, we  fix a constant $r_0\ge 0$ and a natural number $d\ge 1$. Let
$$\C= \big\{\xi=(\xi^1,\cdots,\xi^d): [-r_0,0]\to \R^d \ \text{ cadlag}\big\}.$$
Recall that a path is called cadlag if it is right-continuous
having finite left limits.  We introduce here three different topologies on the space $\C$.
For any $\xi\in\C$, we have
$$\|\xi\|_\infty :=\sum_{i=1}^d\sup_{s\in [-r_0,0]} |\xi^i(s)|<\infty.$$ Under the uniform norm $\|\cdot\|_\infty$, the space $\C$
is complete but not separable. To make $\C$ a Polish space one takes
 the Skorohod metric which is weaker than the uniform metric. Finally, we will also use the topology of  pointwise convergence under  which the space $\C$ is not complete but separable.  The topology of  pointwise convergence is weaker than the other two.   Although these three topologies are different, they all induce the   product $\si$-field on $\C$: $\scr B(\C):= \si\{\gg\mapsto\gg(\theta): \theta\in [-r_0,0]\}$.

For any cadlag $f: [-r_0,\infty)\to \R^d$  and
$t\ge 0$,  we let  $f_t\in\C$ be such that $f_t(\theta)=f(\theta+t)$ for $\theta\in
[-r_0,0]$. We call $(f_t)_{t\ge 0}$ the segment of $(f(t))_{t\ge -r_0}.$
 Next, define $f_{t-}\in\C $ for $t>0$ such that  $f_{t-}(\theta)=f(t+\theta)$ for $\theta\in [-r_0,0)$ and $f_{t-}(0)=f(t-):=\lim_{s\uparrow t} f(s)$. Note that if we replace the definition of $f_{t-}$
 by the usual $\tt f_{t-}(\theta):= f((t+\theta)-)$ for all $\theta\in [-r_0,0]$ then $\tt f_{t-}$ is not necessarily in  $\C$. Moreover, since $\tt f_{t-}$ is a function of $f_{t-}$,
 the equation \eqref{1.0} we consider below also covers the corresponding one for $\tt f_{t-}$ in place of $f_{t-}$.

Now, let $B(t)$ be an $m$-dimensional Brownain motion and let $N(\d s,\d
z)$ be a Poisson counting process with characteristic measure $\nu$
on a measurable space $(E,\scr E)$ with respect to a complete
filtered probability space $(\OO,\F,\{\F_{t}\}_{t\ge 0},\P)$.  We
assume that $B$ and $N$ are independent. We will consider the
order-preservation of SFDEs driven by $B$ and $N$. To characterize
the non-Lipschitz regularity of coefficients in the SFDEs, we
introduce the following class of control functions:
\beg{equation*}\beg{split} \U= \Big\{u\in C^1((0,\infty);[1,\infty)):\ &\int_0^1\ff{\d s}{su(s)} =\infty,\\
& s\mapsto su(s)\  \text{is\ increasing\ and \ concave}\Big\}. \end{split}\end{equation*}
Typical elements in this class are  $u(s)=1$ and $u(s) = \log (\e+s^{-1}).$

Consider the following SFDEs on $\R^d$:
\beq\label{1.0} \beg{cases} \d X(t)= b(t,X_t)\d t+\si(t,X_t)\d B(t)+\int_E \gg(t,X_{t-},z)N(\d t,\d z),\\
\d \bar X(t)= \bar b(t,\bar X_t)\d t+\si(t,\bar X_t)\d B(t)+\int_E \bar \gg(t,\bar X_{t-},z)N(\d t,\d z),\end{cases}\end{equation}where
\beg{equation*}\beg{split} &b,\bar b: [0,\infty)\times \C\times\OO\to \R^d,\ \ \si,\bar\si: [0,\infty)\times\C\times\OO\to \R^d\otimes\R^m,\\
&\gg,\bar\gg: [0,\infty)\times \C\times E\times\OO\to \R^d\end{split}\end{equation*} are  progressively measurable.

For any $s\ge0$ and $\F_s$-measurable $\C$-valued random variables $\xi,\bar\xi$, a solution to (\ref{1.0}) for $t\ge s$ with $(X_s,\bar X_s)= (\xi,\bar\xi)$ is a  cadlag adapted process $(X(t),\bar X(t))_{t\ge s}$  such that $\P$-a.s. for all $t\ge s,$
\beg{equation*}\beg{split} &X(t) = \xi(0)+ \int_s^t  b(r,X_r)\d r+ \int_s^t \si(r, X_r)\d B(r )+\int_{[s,t]\times E} \gg(r,X_{r-},z)N(\d r,\d z),\\
&\bar X(t) = \bar\xi(0)+ \int_s^t \bar b(r,\bar X_r)\d r+ \int_s^t\bar \si(r,\bar X_r)\d B(r) +\int_{[s,t]\times E}\bar \gg(r,\bar X_{r-},z)N(\d r,\d z),\end{split}\end{equation*}where, according to the initial condition $(X_s, \bar X_s)=(\xi,\bar\xi)$, $X_r$ and $\bar X_r$ for $r\ge s$ are well defined.

To ensure the existence and uniqueness of solutions, we make use of the following assumptions:    \beg{enumerate} \item[{\bf (A1)}] There exist some positive $K\in C([0,\infty))$ and $u\in\U$ such that $\P$-a.s.
 \beg{equation*}\beg{split} &|b(t,\xi)-b(t,\eta)|^2+|\bar b(t,\xi)-\bar b(t,\eta)|^2+\|\si(t,\xi)-\si(t,\eta)\|_{HS}^2 + \|\bar\si(t,\xi)-\bar\si(t,\eta)\|_{HS}^2\\
&\qquad + \int_E \big(|\gg(t,\xi,z)-\gg (t,\eta,z)|^2+ |\bar \gg(t,\xi,z)-\bar \gg (t,\eta,z)|^2\big)\nu(\d z)\\
&\qquad +\bigg(\int_E \big(|\gg(t,\xi,z)-\gg (t, \eta,z)|+  |\bar \gg(t,\xi,z)-\bar \gg (t,\eta,z)|\big)\nu(\d z)\bigg)^2\\
&\le K(t)   \|\xi -\eta \|_\infty^2 u(\|\xi -\eta \|_\infty^2),\ \ \xi, \eta \in\C, t\ge 0,\end{split}\end{equation*}where $\|\cdot\|_{HS}$ denotes the Hilbert-Schmidt norm.
 \item[{\bf (A2)}] For any $T>0$ there exists a constant $C(T)>0$ such that $\P$-a.s.
  \beg{equation*}\beg{split} &\sup_{t\in [0,T]} \big(|b(t,0)|^2+|\bar b(t,0)|^2 +\|\si(t,0)\|_{HS}^2+\|\bar\si(t,0)\|_{HS}^2 \big)\\
  &+\int_{[0,T]\times E}\big( |\gg(t,0,z)|^2+|\bar\gg(t,0,z)|^2\big)\d t\nu(\d z)\le C(T).\end{split}\end{equation*}
 \end{enumerate}
If $u\equiv 1$ then {\bf (A1)} reduces to the usual Lipschitz condition. In general,  {\bf (A1)} allows the coefficients to be non-Lipschitzian.

According to Theorem \ref{T3.1} below,   for any $s\ge 0$ and $\F_s$-measurable $\C$-valued random variables $\xi,\bar\xi$,     the equation  (\ref{1.0}) has a unique solution  for $t\ge s$ with $X_s=\xi$ and $\bar X_s=\bar\xi$ and the solution is non-explosive. We denote the solution  by   $\{X(s,\xi;t), \bar X(s,\bar \xi;t)\}_{t\ge s}. $

To introduce the notion of order-preservation of the solutions, we take the usual partial-order on $\R^d$; i.e., for $x=(x^1,\cdots, x^d),y=(y^1,\cdots, y^d)\in\R^d$, we write $x\le y$ if $x^i\le y^i$ holds for all $1\le i\le d.$
Similarly,  for $\xi=(\xi^1,\cdots,\xi^d),\eta=(\eta^1,\cdots,\eta^d)\in\C$, we write $\xi\le \eta$ if $\xi^i(\theta)\le \eta^i(\theta)$ holds for all $\theta\in [-r_0,0]$ and $1\le i\le d.$
Moreover, for any $\xi,\eta\in\C$, let $\xi\land\eta\in \C$ be such that
$(\xi\land\eta)^i=\min\{\xi^i,\eta^i\}, 1\le i\le d$ and let
$\xi\lor \eta= -\{(-\xi)\land(-\eta)\}.$

\beg{defn} The solutions of $(\ref{1.0})$ are called order-preserving if, for any $s\ge 0$ and $\F_s$-measurable $\C$-valued random variables $\xi,\bar\xi$  with $\P$-a.s.  $\xi\le\bar\xi$,   $\P$-a.s. for all $t\ge s$ it holds that $X(s,\xi;t)\le \bar X(s,\bar\xi;t)$.   \end{defn}

\beg{thm}\label{T1.1} Assume {\bf (A1)} and {\bf (A2)}. The solutions to $(\ref{1.0})$ are order-preserving if the following   conditions are satisfied:
\beg{enumerate} \item[$(1)$]  $\P$-a.s. for all   $1\le i\le d, t\ge 0$ and  $\xi,\bar\xi\in\C$ with  $\xi\le\bar\xi$ and $\xi^i(0)=\bar\xi^i(0),$ it holds that $b^i(t,\xi)\le \bar b^i(t,\bar\xi)$.
\item[$(2)$]  $\P$-a.s.  for all   $t\ge 0, 1\le i\le d, 1\le j\le m$ and $\xi,\bar\xi\in \C$ with  $\xi^i(0)=\bar\xi^i(0)$, it holds  that $\si^{ij}(t,\xi)=\bar\si^{ij}(t,\bar \xi)$.
    \item[$(3)$]   $\nu\times\P$-a.e. for all all $1\le i\le d, t\ge 0$ and $\xi,\bar\xi\in\C$ with  $\xi\le \bar\xi$, it holds  that  $\xi^i(0)+\gg^i(t,\xi,\cdot)\le \bar\xi^i(0)+\bar\gg^i(t,\bar\xi,\cdot)$. \end{enumerate}
\end{thm}

Note that condition (2) implies that, $\P$-a.s. for every $1\le i\le d, 1\le j\le m$ and $t\ge 0$,  $\si^{ij}(t,\xi)$ depends only on $\xi^i(0)$, so that   the diffusion coefficient does not contain any delay.
Comparing with existing comparison theorems derived in the above mentioned references for one-dimensional equations, Theorem \ref{T1.1} has a  rather broad range of applications. Next, a multidimensional comparison theorem without delay has been presented in \cite[Theorem 296]{ST} where the  condition $2^0$ implies that $\bar b^i(t,x)$ depends only on $x^i$, and is thus much stronger than the condition  (1) in Theorem \ref{T1.1} with $r_0=0$ (i.e. the case without delay).
Moreover, when $r_0=0$ (i.e. without delay) Theorem \ref{T1.1} also covers the comparison theorem derived recently in \cite{Z} for $\nu(E)<\infty$ and Lipschitzian coefficients. Finally, when $b$ and $\bar b$ are deterministic without delay, condition (1) coincides with condition $(C_{a,b})$ in \cite{Ass}.

On the other hand, our next result    shows that the conditions  in Theorem \ref{T1.1} are also necessary for the order-preservation  under mild assumptions. This result is new even in the case of coefficients without delay.

\beg{thm}\label{T1.2}  Let $\C_\infty$ and $\C_p$ denote the space $\C$ equipped with the uniform metric and the topology of pointwise convergence, respectively. Assume  {\bf (A1)},  {\bf (A2)} and that the solutions to $(\ref{1.0})$ are order-preserving.  \beg{enumerate}\item[{\rm (I)}] If
$\P$-a.s. $b,\bar b\in C([0,\infty)\times \C_p; \R^d)$ and,  for any $n\ge 1$,
\beq\label{AB} \lim_{\vv\downarrow 0} \sup_{t\in [0,n], \|\xi\|_\infty \le n} \int_{E}  \vv\land (|\gg(t,\xi,z)|+|\bar\gg(t,\xi,z)|)\nu(\d z) =0\end{equation} then condition $(1)$ holds.
\item[{\rm (II)}] If $\P$-a.s.
$\si ,\bar\si \in C([0,\infty)\times\C_\infty; \R^d\otimes\R^m)$ and  $\nu\times \P$-a.e. $\gg,\bar\gg\in C([0,\infty)\times\C_p; \R^d)$  then  condition $(2)$ holds.
\item[{\rm (III)}] If $\nu\times\P$-a.e. $\gg,\bar\gg\in C([0,\infty)\times\C_p; \R^d)$   then  condition $(3)$ holds. \end{enumerate}    \end{thm}
Note that condition (\ref{AB}) holds if either $\nu$ is finite or $\gg(t,\xi,\cdot),\bar\gg(t,\xi,\cdot)$ are integrable with respect to $\nu$ locally uniformly in $(t,\xi)$.

In the next section we present proofs of the above two theorems. In Section 3, we present a  result  on the existence and uniqueness of solutions to stochastic functional equations with jumps.

\section{Proofs of Theorems \ref{T1.1} and \ref{T1.2}}

\beg{proof}[Proof of Theorem \ref{T1.1}] Assume that (1)--(3) hold. For any $t_0\ge 0$ and $\F_{t_0}$-measurable $\C$-valued random variables $\xi,\bar\xi$  such that $\P$-a.s. $\xi\le\bar\xi$, we aim to prove
that for any $T>t_0$, \beq\label{P} \E  \sup_{r\in [t_0,T]}(X^i(t_0,\xi; r)-\bar X^i(t_0,\bar\xi;r))^+=0,\ \ 1\le i\le d.\end{equation} For equations without delay, this can be done by using a Tanaka type formula (see \cite[153]{ST}). Below, we shall adopt an approximation argument which plays the same role as the Tanaka type formula. For simplicity, we will  denote $X(t)=X(t_0,\xi;t)$ and $\bar X(t)=\bar X(t_0,\bar\xi;t)$ for $t\ge t_0-r_0$. Recall that $X(t_0,\xi;t)=\xi(t-t_0)$ and $\bar X(t_0,\bar\xi;t)= \bar \xi(t-t_0)$ for $t\in [t_0-r_0,t_0].$

For any $n\ge 1,$ let $\psi_n: \R\to [0,\infty)$ be constructed as follows: $\psi_n(s)=\psi_n'(s)=0$ for $s\in (-\infty,0]$, and
$$\psi_n''(s)=\beg{cases} 4n^2s, & s\in [0,\ff 1 {2n}],\\
-4n^2(s-\ff 1 n), & s\in [\ff 1 {2n}, \ff 1 n],\\
0, &\text{otherwise}.\end{cases}$$  We have
\beq\label{1.3} 0\le \psi_n'\le 1_{(0,\infty)}, \ \text{and\ as\ } n\uparrow\infty: \ 0\le \psi_n(s)\uparrow s^+,\ \ s\psi_n''(s)\le  1_{(0,\ff 1 n)}(s)\downarrow 0.\end{equation}
Let $$\tau_k=\inf\{t\ge t_0: |X(t)-X(t)\land\bar X(t)|\ge k\},\ \ k\ge 1.$$ Since $\si=\bar\si$  by (2) and since $$\psi_n(X^i(t_0)-\bar X^i(t_0))=\psi_n(\xi^i(0)-\bar\xi^i(0))=0,$$   by $\xi\le \bar\xi,$  It\^o's formula yields $\P$-a.s.
\beq\label{1.4} \beg{split}&\psi_n(X^i(t\land \tau_k)-\bar X^i(t\land\tau_k))^2\\
&= M_i(t\land\tau_k)+2\int_{t_0}^{t\land\tau_k} (b^i(s,X_s)-\bar b^i(s,\bar X_s))\{\psi_n\psi_n'\}(X^i(s)-\bar X^i(s))\d s\\
&\quad +   \sum_{j=1}^m \int_{t_0}^{t\land\tau_k} \big(\si^{ij}(s,X_s)- \si^{ij}(s,\bar X_s)\big)^2 \{\psi_n\psi_n''+{\psi_n'}^2\}(X^i(s)-\bar X^i(s))\d s\\
&\quad + \int_{[t_0,t\land\tau_k]\times E}\Big\{\psi_n\big(X^i(s-)-\bar X^i(s-)+ \gg^i(s,X_{s-},z)- \bar\gg^i(s,\bar X_{s-},z)\big)^2\\
&\qquad\qquad\qquad\qquad -\psi_n\big(X^i(s-)-\bar X^i(s-)\big)^2\Big\}N(\d s,\d z) \end{split}\end{equation}for any $k,n\ge 1$, $1\le i\le d$ and $t\ge t_0$
where
$$M_i(t):= 2\sum_{j=1}^m \int_{t_0}^t \big(\si^{ij}(s, X_s)- \si^{ij}(s,\bar X_s)\big)\{\psi_n\psi_n'\}(X^i(s)-\bar X^i(s))\d B^j(s).$$
Noting that $0\le\psi_n'(X^i(s)-\bar X^i(s))\le 1_{\{X^i(s)>\bar X^i(s)\}}$ and when $X^i(s)>\bar X^i(s)$ one has
$(X_s\land \bar X_s)^i(0)=(\bar X_s)^i(0),$ it follows from (1) that $\P$-a.s.
$$(b^i(s,X_s\land\bar X_s)-\bar b^i(s,\bar X_s))\{\psi_n\psi_n'\}(X^i(s)-\bar X^i(s))\le 0,\ \ n\ge 1, s\in [t_0,T].$$ Combining this with
 {\bf (A1)} and $0\le\psi_n'\le 1$, we obtain $\P$-a.s.
\beq\label{1.5} \beg{split} &2\int_{t_0}^{t\land \tau_k}(b^i(s,X_s)-\bar b^i(s,\bar X_s))\{\psi_n\psi_n'\}(X^i(s)-\bar X^i(s))\d s\\
&=2\int_{t_0}^{t\land \tau_k}\Big[ (b^i(s,X_s)- b^i(s,X_s\land \bar X_s))\{\psi_n\psi_n'\}(X^i(s)-\bar X^i(s))\\
&\qquad\qquad + (b^i(s,X_s\land\bar X_s)-\bar b^i(s,\bar X_s))\{\psi_n\psi_n'\}(X^i(s)-\bar X^i(s))\Big]\d s\\
&\le 2\int_{t_0}^{t\land \tau_k}|b^i(s,X_s)- b^i(s, X_s\land \bar X_s)|\cdot \psi_n(X^i(s)- \bar X^i(s))\d s\\
&\le \int_{t_0}^{t\land \tau_k}\Big[8T|b^i(s,X_s)- b^i(s, X_s\land \bar X_s)|^2+\ff 1 {8T} \psi_n(X^i(s)-\bar X^i(s))^2\Big]\d s \\
&\le C(T) \int_{t_0}^{t\land \tau_k}\|X_s-X_s\land \bar X_s\|_\infty^2 u(\|X_s-X_s\land\bar X_s\|_\infty^2)\d s \\
&\quad +\ff 1 {8} \sup_{s\in [t_0,t\land \tau_k]}\psi_n(X^i(s)-\bar X_i(s))^2,\ \ \ t\in [t_0,T],\ n,k\ge 1 \end{split}\end{equation} for some constant $C(T)>0$. Next, since (2) implies that $\si^{ij}(s, X_s)=\bar\si^{ij}(s, X_s)$ depends only on $X^i(s)$, it follows from  {\bf (A1)} and  (\ref{1.3})  that  $\P$-a.s.
\beq\label{1.6} \beg{split}&\sum_{j=1}^m |\si^{ij}(s,X_s)-\bar \si^{ij}(s,\bar X_s)|^2\{\psi_n\psi_n''+{\psi_n'}^2\}X^i(s)-\bar X^i(s))\\
&\le C(T) 1_{\{X^i(s)-\bar X^i(s)\in (0,\ff 1 n)\}} |X^i(s)-\bar X^i(s)|^2 u(|X^i(s)-\bar X^i(s)|^2) \\
&\qquad + C(T) \|X_s-X_s\land \bar X_s\|_\infty^2 u(\|X_s-X_s\land\bar X_s\|_\infty^2)\\
&\le  \vv(n)+ C(T) \|X_s-X_s\land \bar X_s\|_\infty^2 u(\|X_s-X_s\land\bar X_s\|_\infty^2),\ \ n\ge 1, s\in [t_0,T]\end{split}\end{equation} for some constant $C(T)>0$ where, since $u\in\scr U$,
$$\vv(n):= C(T)\sup_{s\in (0,n^{-2})} s u(s) \downarrow 0\ \text{as}\ n\uparrow\infty.$$
Moreover, {\bf (A1)}, (\ref{1.3})  and (2) also imply $\P$-a.s.
\beq\label{1.7}  \beg{split} &\sum_{j=1}^m  |\si^{ij} (s,X_s)-\bar\si^{ij}(s,\bar X_s)|^2\{\psi_n\psi_n'\}^2(X^i(s)-\bar X^i(s)) \\
&\le C_1(T) \|X_s-X_s\land \bar X_s\|_\infty^2 u(\|X_s-X_s\land\bar X_s\|_\infty^2)\psi_n(X^i(s)-\bar X^i(s))^2\end{split}\end{equation} for some constant $C_1(T)>0$ and all $ n\ge
1, s\in [t_0,T]$  so that, by Burkholder-Davis-Gundy's inequality,
\beq\label{QQ} \beg{split} &\E\sup_{s\in [t_0,t]}  M_i(s\land \tau_k)   
 \le C_2(T) \\
 &\times \E \bigg(\int_{t_0}^{t\land\tau_k} \|X_s-X_s\land \bar X_s\|_\infty^2 u(\|X_s-X_s\land\bar X_s\|_\infty^2)\psi_n(X^i(s)-\bar X^i(s))^2\d s\bigg)^{\ff 1 2}\\
 &\le C_3(T)  \E \int_{t_0}^{t\land\tau_k} \|X_s-X_s\land \bar X_s\|_\infty^2 u(\|X_s-X_s\land\bar X_s\|_\infty^2) \d s\\
&\qquad + \ff 1 8 \E\sup_{s\in [t_0, t\land\tau_k]}\psi_n(X^i(s)-\bar X_i(s))^2 \end{split}
\end{equation} holds for some  constants $C_2(T),C_3(T)>0$ and all $n\ge 1, t\in [t_0,T].$

Finally, by (3), we have
\beq\label{A**B}X^i(s)\land\bar X^i(s)+\gg^i(s, X_{s}\land\bar X_{s}, \cdot)\le \bar X^i(s)+\bar\gg^i (s, \bar X_{s}, \cdot),\ \  \nu\times\P\text{-a.e.}\end{equation} If $X^i(s)\le \bar X^i(s)$, then (\ref{A**B}) becomes
$$X^i(s)+\gg^i(s, X_s\land\bar X_s,\cdot)\le \bar X^i(s)+\bar \gg^i(s, \bar X_s,\cdot),\ \ \nu\times\P\text{-a.e.},$$ so that $0\le\psi_n'\le 1$ and $\psi_n(s)=0$ for $s\le 0$  imply
\beg{equation*}\beg{split} &\psi_n\big(X^i(s)-\bar X^i(s)+
\gg^i(s,X_{s},\cdot)- \bar\gg^i(s,\bar X_{s},\cdot)\big)^2
-\psi_n(X^i(s)-\bar X^i(s))^2\\
&=\psi_n\big(X^i(s)-\bar X^i(s)+ \gg^i(s,X_{s},\cdot)- \bar\gg^i(s,\bar
X_{s},\cdot)\big)^2 \\
&=\psi_n\big( \gg^i(s,X_{s},\cdot)-\gg^i(s, X_{s}\land\bar X_{s}, \cdot)\\
&\qquad \quad+(X^i(s)+\gg^i(s, X_{s}\land\bar X_{s}, \cdot))- (\bar
X^i(s)+\bar\gg^i(s,\bar X_{s},\cdot))\big)^2 \\
&\le \psi_n\big( \gg^i(s,X_{s},\cdot)-\gg^i(s, X_{s}\land\bar X_{s},
\cdot)\big)^2\\
&\le|\gg^i(s,X_{s},\cdot)-\gg^i(s, X_{s}\land\bar X_{s},
\cdot)|^2,\ \ \nu\times\P\text{-a.e.}\end{split}\end{equation*} On the other hand,  if $X^i(s)\ge \bar
X^i(s)$ then (\ref{A**B}) becomes
$$\gg^i(s,X_s\land\bar X_s,\cdot)\le \bar \gg^i(s,\bar X_s,\cdot),\ \ \nu\times\P\text{-a.e.},$$
so that $0\le\psi_n'\le 1$ implies that $\nu\times\P\text{-a.e.}$
  \beg{equation*}\beg{split}
&\psi_n\big(X^i(s)-\bar X^i(s)+ \gg^i(s,X_{s},\cdot)- \bar\gg^i(s,\bar
X_{s},\cdot)\big)^2
-\psi_n(X^i(s)-\bar X^i(s))^2\\
&=\psi_n\big(X^i(s)+ \gg^i(s,X_{s},\cdot)-\bar
X^i(s)-\gg^i(s, X_{s}\land\bar X_{s}, \cdot) \\
&\quad+ \gg^i(s, X_{s}\land\bar X_{s}, \cdot) -  \bar\gg^i(s,\bar X_{s},\cdot) \big)^2-\psi_n(X^i(s)-\bar X^i(s))^2 \\
&\le \psi_n\big(X^i(s)-\bar X^i(s)+ \gg^i(s,X_{s},\cdot)-\gg^i(s,
X_{s}\land\bar X_{s},
\cdot))^2-\psi_n(X^i(s)-\bar X^i(s))^2\\
&\le 2|\gg^i(s,X_{s},\cdot)-\gg^i(s, X_{s}\land\bar X_{s},
\cdot)|\\
&\qquad \times \{\psi_n(X^i(s)-\bar X^i(s))+ |\gg^i(s,X_{s},\cdot)-\gg^i(s, X_{s}\land\bar X_{s},
\cdot)|\}.\end{split}\end{equation*}
In conclusion, we have
\beg{equation*}\beg{split}
&\psi_n\big(X^i(s)-\bar X^i(s)+ \gg^i(s,X_{s},\cdot)- \bar\gg^i(s,\bar
X_{s},\cdot)\big)^2
-\psi_n(X^i(s)-\bar X^i(s))^2\\
&\le 2|\gg^i(s,X_{s},\cdot)-\gg^i(s, X_{s}\land\bar X_{s},
\cdot)|^2
  +  2 \psi_n(X^i(s)-\bar X^i(s))   |\gg^i(s,X_{s},\cdot)-\gg^i(s, X_{s}\land\bar X_{s},
\cdot)| .\end{split}\end{equation*}
Combining this with {\bf (A1)} we obtain
 \beq\label{1.7'} \beg{split} &\E\sup_{r\in [t_0,t]} \int_{[t_0,r\land\tau_k]\times E} \Big\{\psi_n\big(X^i(s-)-\bar X^i(s-)+ \gg^i(s,X_{s-},z)- \bar\gg^i(s,\bar X_{s-},z)\big)^2\\
 &\qquad\qquad\qquad \qquad-\psi_n(X^i(s-)-\bar X^i(s-))^2\Big\}^+  N(\d s,\d z)\\
 &= \E   \int_{[t_0,t\land\tau_k]\times E} \Big\{\psi_n\big(X^i(s-)-\bar X^i(s-)+ \gg^i(s,X_{s-},z)- \bar\gg^i(s,\bar X_{s-},z)\big)^2\\
 &\qquad\qquad\qquad \qquad-\psi_n(X^i(s-)-\bar X^i(s-))^2\Big\}^+  N(\d s,\d z)\\
&= \E\int_{[t_0,t\land\tau_k]\times E} \Big\{\psi_n\big(X^i(s)-\bar X^i(s)+ \gg^i(s,X_{s},z)- \bar\gg^i(s,\bar X_{s},z)\big)^2\\
 &\qquad\qquad\qquad \qquad-\psi_n(X^i(s)-\bar X^i(s))^2\Big\}^+  \d s\nu(\d z)\\
  &\le C(T) \E \int_{t_0}^{t\land\tau_k}\|X_{s}-X_{s}\land\bar X_{s}\|_\infty ^2u(\|X_{s}-X_{s}\land\bar X_{s}\|_\infty^2) \d s\\
&\qquad + \ff 1 4 \E \sup_{s\in [t_0,t\land\tau_k]} \psi_n(X^i(s)-\bar X^i(s))^2\end{split}\end{equation} for some constant $C(T)>0$ and all $ n\ge 1, t\in [t_0,T].$

Now, let
$$\phi_k(s)=\sup_{r\in [t_0-r_0, s\land\tau_k]}|X(r)-X(r)\land \bar X(r)|^2,\ \ s\ge t_0.$$
 By substituting  \eqref{1.5}, \eqref{1.6}, \eqref{QQ}, \eqref{1.7'} into (\ref{1.4}), and noting that  $X_{t_0}\le\bar X_{t_0}$, we obtain 
 \beg{equation*}\beg{split} &\E  \sup_{r\in [t_0-r_0, t\land \tau_k]} \psi_n(X^i(r)-\bar X^i(r))^2=   \E  \sup_{r\in [t_0, t\land \tau_k]} \psi_n(X^i(r)-\bar X^i(r))^2\\
&\le C(T)  \int_{t_0}^t \E \big\{\phi_k(s)u(\phi_k(s))\big\}\d s +\vv(n) +\ff 1 2 \E  \sup_{r\in [t_0-r_0, t\land \tau_k]} \psi_n(X^i(r)-\bar X^i(r))^2\end{split}\end{equation*} 
for some constants $C(T)>0$ and $\vv(n)>0$ with $\vv(n)\to 0$ as $n\to\infty,$ and all $k,n\ge 1, t\in [t_0,T], 1\le i\le d.$ Therefore, for any $n,k\ge 1$ and $t\in [t_0,T],$ we have 
$$ \sum_{i=1}^d\E  \sup_{r\in [t_0-r_0, t\land \tau_k]} \psi_n(X^i(r)-\bar X^i(r))^2 \le 2 d  C(T)  \int_{t_0}^t \E \big\{\phi_k(s)u(\phi_k(s))\big\}\d s
+ 2 d \vv(n).$$
Letting  $n\uparrow\infty$ and using  Jensen's inequality, we arrive at
$$\E  \phi_k(t)\le 2 d C(T) \int_{t_0}^t \{ \E\phi_k(s)\} u\big(\E \phi_k(s)\big)\d s,\ \ t\in [t_0,T], k\ge 1.$$ Since $\int_0^1\ff 1 {su(s)}\d s=\infty$, by   Bihari's inequality this implies that (see e.g. the end of the proof of Theorem 4.2 in \cite{SWY})
$$ \E \phi_k(T)=0,\ \ k\ge 1.$$Letting $k\uparrow\infty$   proves (\ref{P}). \end{proof}

To prove Theorem \ref{T1.2}, we need the following Lemma \ref{L2.1}. For any $h\in C^2_b(\R^d),$ let
\beg{equation*}\beg{split} &(Lh)(t,\xi)=\sum_{i=1}^d b^i(t,\xi)\pp_i h(\xi(0)) +\ff 1 2 \sum_{i,j=1}^d(\si\si^*)(t,\xi) \pp_i\pp_j h(\xi(0))\\
&\qquad\qquad\qquad + \int_{E} \big\{h(\xi(0)+\gg(t,\xi,z))-h(\xi(0))\big\}\nu(\d z),\\
&(\bar Lh)(t,\xi)=\sum_{i=1}^d \bar b^i(t,\xi)\pp_i h(\xi(0)) +\ff 1 2 \sum_{i,j=1}^d(\bar \si\bar\si^*)(t,\xi) \pp_i\pp_j h(\xi(0))\\
&\qquad\qquad\qquad + \int_{E}
\big\{h(\xi(0)+\bar\gg(t,\xi,z))-h(\xi(0))\big\}\nu(\d z), \ t\ge 0,
\xi\in\C,\end{split}\end{equation*} where $\pp_i (1\le i\le d)$ is
the derivative with respect to the $i$-th component in $\R^d$. By {\bf (A1)} and {\bf (A2)}, $Lh$ and $\bar Lh$ are locally bounded with respect to the usual metric on $[0,\infty)$ and the uniform norm on $\C$.

Let
$\M$ be the class of all   increasing functions on $\R^d$ where a
function $h$ on $\R^d$ is called increasing if $h(x)\le h(y)$
holds for all $x\le y.$

\beg{lem}\label{L2.1} Assume {\bf (A1)} and {\bf (A2)}.   If the solutions to $(\ref{1.0})$ are order-preserving then, for any $s\ge 0$, $\xi,\bar\xi\in\C$  with   $\xi\le\bar\xi$,  and   $h\in \M\cap C^2_b(\R^d)$ with $h(\xi(0))=h(\bar \xi(0)),$ it holds that $\P$-a.s.
$$\E\Big(\liminf_{t\downarrow s} (Lh)(t,X_t(s,\xi))\Big|\F_s\Big)\le \E\Big(\limsup_{t\downarrow s} (\bar Lh)(t,\bar X_t(s,\bar\xi))\Big|\F_s\Big)$$
where $X_\cdot(s,\xi)$ and $\bar X_\cdot(s,\bar\xi)$ are the segment processes of $X(s,\xi;\cdot)$ and $\bar X(s;\bar\xi,\cdot)$, respectively. \end{lem}

\beg{proof}   Simply denote $X(t)= X(s,\xi;t),\bar X(t)=\bar
X(s,\bar\xi;t)$ for $t\ge s-r_0$. Let $$\tau= \inf\{t\ge s: |X(t)|+|\bar
X(t)|\ge 1+\|\xi\|_\infty+\|\bar\xi\|_\infty\}.$$
 Since $h(\xi(0))=h(\bar\xi(0))$ and $X(t)\le\bar X(t)$ for all $t\ge s,$ we have
\beq\label{*} \E \big(h(X(t\land\tau))\big|\F_s\big)-h(\xi(0)) \le \E\big( h(\bar X
(t\land\tau))\big|\F_s\big)-h(\bar\xi(0)),\ \ t\ge s.\end{equation}
  By  It\^o's formula and   Fatou's lemma, it is easy to see that
\beg{equation}\label{*D}\beg{split} &\liminf_{t\downarrow s} \ff{\E \big(h(X(t\land\tau))\big|\F_s\big)- h(\xi(0))}{t-s} \ge \E\Big(\liminf_{t\downarrow s}(L h) (t,X_t)\Big|\F_s\Big),\\
&\limsup_{t\downarrow s} \ff{\E \big(h(\bar X(t\land\tau))\big|\F_s\big)- h(\bar \xi(0))}{t-s} \le
\E\Big(\limsup_{t\downarrow s}  (\bar L h) (t,\bar X_t)\Big|\F_s\Big).\end{split}\end{equation} Combining this with (\ref{*}) we finish the proof.

Below we only prove the first formula in (\ref{*D}), as the proof of the second is completely similar.  By   It\^o's formula, for $t\in (s,s+1],$ we have
\beg{equation*}\beg{split} &\ff{\E\big(h(X(t\land\tau))\big|\F_s\big) - h(\xi(0))}{ t-s} = \ff{\E\big(\int_s^{t\land\tau}
 (Lh)(r,X_r)\big|\F_s\big)}{t-s}\\
&\ge \E\Big(1_{\{\tau>t\}}\inf_{r\in (s,t\land \tau]}(Lh)(r, X_r)\Big|\F_s\Big)-C\P(\tau<t|\F_s)\end{split}\end{equation*}where, due to {\bf (A1)} and {\bf (A2)},
$$C:= \sup_\OO \big\{|Lh|(r,\eta):\ r\in [s,s+1], \|\eta\|_\infty\le 1+\|\xi\|_\infty+\|\bar\xi\|_\infty\big\}<\infty.$$
Since  $\tau>s$ due to the right-continuity of the solution  and since  $Lh$ is locally bounded, by letting $t\downarrow s,$ we obtain the first formula in (\ref{*D}) from  Fatou's lemma.
\end{proof}

\beg{proof}[Proof of Theorem \ref{T1.2}]  Let $1\le i\le d$  and $t_0\ge 0$ be fixed. For any  $\xi,\bar\xi\in\C$, let  $X(t)=X(t_0,\xi,;t), \bar X(t)= \bar X(t_0,\bar\xi;t).$

(a) Proof of   (III). Let     $\xi\le \bar\xi$. {\bf (A1)} and {\bf (A2)} imply that $b,\si$ and
$\int_{E}\big(|\gg(\cdot,\cdot,z)|+|\bar\gg(\cdot,\cdot,z)|\big)\nu(\d z)$ are locally bounded. We aim to
prove \beq\label{C*D} \xi^i(0)+\gg^i(t_0,\xi,\cdot)\le
\bar\xi^i(0)+\bar\gg^i(t_0,\bar\xi,\cdot), \,\,\nu\times\P\text{-a.e.} \end{equation} Due to the continuity of
$\gg$ and $\bar\gg$ in the first two variables and the separability of $[0,\infty)\times\C_p$,   this implies condition   (3).

Let $\tau=\inf\{t\ge t_0: \|X_t-\xi\|_\infty+\|\bar X_t-\bar\xi\|_\infty \ge 1\}$ and let $\psi_n$ be as in the proof of Theorem \ref{T1.1}.
By
It\^o's formula and the local boundedness of the coefficients and the
right-continuity of $X(s)$,  for any $t>t_0,$ we have
\beg{equation*}\beg{split} &\E \psi_n\big(X^i(t\land\tau)-\bar X^i(t\land\tau)\big)\\
&= \E\int_{t_0}^{t\land\tau} \Big\{(b^i(s,X_s)-\bar b^i(s,\bar X_s))\psi_n'\big(X^i(s)-\bar X^i(s)\big)\\
 &\quad +\ff 1 2 \sum_{j=1}^m (\si^{ij}(s,X_s)-
\bar\si^{ij}(s,\bar X_s))^2\psi_n''(X^i(s)-\bar X^i(s)) \\
&\quad +  \int_{E}\big[\psi_n\big(X^i(s)-\bar X^i(s)+
\gg^i(s,X_{s},z)- \bar\gg^i(s,\bar X_{s},z)\big)
-\psi_n(X^i(s)-\bar X^i(s))\big]\nu(\d z)\Big\}\d s.\\
\end{split}\end{equation*} Since
$X^i(s)\le\bar X^i(s), X^i(t\land\tau)\le\bar X^i(t\land\tau)$ and $\psi_n(s)=\psi_n'(s)=\psi_n''(s)=0$ for $s\le 0$, this implies
 \beg{equation*}\beg{split}
\E\int_{t_0}^{t\land\tau}\Big\{\int_{E}\psi_n\big(X^i(s)-\bar
X^i(s)+ \gg^i(s,X_{s},z)- \bar\gg^i(s,\bar X_{s},z)\big) \nu(\d
z)\Big\}\d s=0, \, \ t\ge t_{0}.
\end{split}\end{equation*} By $\tau>t_0$,  the right-continuity of the solutions (which implies  $\lim_{s\downarrow 0} (X_s,\bar X_s)=(\xi,\bar\xi)$ in $\C_p\times\C_p$) and the joint-continuity of $\gg^i$ and $\bar\gg^i$ in the first two variables,   from this and   Fatou's lemma we conclude that
 $$\E \int_{E}\psi_n\big(\xi^i(0)-\bar
\xi^i(0)+ \gg^i(t_{0},\xi,z)- \bar\gg^i(t_{0},\bar
\xi,z)\big) \nu(\d z) =0.$$ Since $\psi_n(s)\uparrow s^+$ as $n\uparrow\infty$, by letting $n\uparrow\infty$ we arrive at
$$ \E\int_{E} \big(\xi^i(0)-\bar
\xi^i(0)+ \gg^i(t_{0},\xi,z)- \bar\gg^i(t_{0},\bar
\xi,z)\big)^+ \nu(\d z) =0,$$
and hence (\ref{C*D}) holds.

(b) Proof of (I).   Let    $\xi\le \bar\xi$ and $\xi^i(0)= \bar\xi^i(0)$. For any $\vv\in (0,1)$, let $\phi_\vv\in C_0^\infty(\R)$ such that $$0\le\phi_\vv\le 1,\ \  \phi_\vv|_{[-\vv,\vv]}=1,\ \
\phi_\vv|_{[-2\vv,2\vv]^c}=0.$$ Take
$$h_\vv(x)= \int_0^{x^i-\xi^i(0)} \phi_\vv (s)\d s,\ \ \ x\in \R^d.$$ Then $h_\vv\in \scr M\cap C_b^2(\R^d)$ and
$$\<y, \nn h_\vv(x)\>:=\sum_{j=1}^d y^j\pp_{x^j} h_\vv(x) = y^i,\ \ x,y\in\R^d, |x^i-\xi^i(0)|\le\vv.$$
So, by the continuity of $b,\bar b$ and the right-continuity of the solutions, we have
\beg{equation*}\beg{split} &\lim_{t\downarrow t_0} \<b(t,X_t),\nn h_\vv(X(t))\>= b^i(t_0,\xi),\ \ \lim_{t\downarrow t_0} \<\bar b(t,\bar X_t),\nn h_\vv(\bar X(t))\>= \bar b^i(t_0,\bar \xi),\\
 &\lim_{t\downarrow t_0} \nn^2 h_\vv(X(t))=\lim_{t\downarrow t_0} \nn^2h_\vv(\bar X(t))=0,\\
& |h_\vv(X(t)+\gg(t,X_t,z))- h_\vv(X(t))| + |h_\vv(\bar X(t)+\bar \gg(t,\bar X_t,z))- h_\vv(\bar X(t))|\\
&\le (4\vv)\land (|\gg(t,X_t,z)|+|\bar\gg(t,\bar X_t,z)|).\end{split}\end{equation*} Combining this with Lemma \ref{L2.1}, we obtain $\P$-a.s.
$$b^i(t_0,\xi)\le \bar b^i(t_0,\bar\xi) +\sup_{t\in [t_0, t_0+1], \|\eta\|_\infty\lor\|\bar\eta\|_\infty\le 1+\|\xi\|_\infty\lor\|\bar\xi\|_\infty}\int_{E}\big\{(4\vv)\land (|\gg(t,\eta,z)|+|\bar\gg(t,\bar\eta,z)|)\big\}\nu(\d z).$$ So, when  $\vv\to 0$,   (\ref{AB}) yields
$b^i(t_0,\xi)\le \bar b^i(t_0,\bar\xi)\ \P$-a.s.  This implies condition (1) by the continuity of $b,\bar b$ and the separability of $[0,\infty)\times\C$.

(c) Proof    (II). If condition (2) does not hold then there exist $1\le i\le d, t_0>0$ and  $\xi,\bar\xi\in\C$ with   $\xi^i(0)= \bar\xi^i(0)$ such that
for some  $ 1\le j\le m$ one has   $ \P(\si^{ij}(t_0,\xi)\ne \bar\si^{ij}(t_0,\bar\xi))>0.$ Since $\si$ and $\bar\si$ are continuous on $[0,\infty)\times\C_\infty,$ there exists a   constant $\vv\in (0,1)$ such that $\P(A_\vv)>0$, where
$$A_{\vv}:=\big\{|\si^{ij}(t,\eta)-\bar\si^{ij}(t,\bar\eta)|\ge \vv \  \text{for}\  t\in [t_0,t_0+\vv], \|\eta-\xi\|_\infty+\|\bar\eta-\bar\xi\|_\infty\le\vv\big\}.$$
Let
\beg{equation*}\beg{split}&\tt \tau=\inf\{t\ge t_0: |\si^{ij}(t,X_t)-\bar\si^{ij}(t,\bar X_t)|\le\vv \}\\
& \tau=\inf\{t\ge t_0:   \|X_t-\xi\|_\infty +\|\bar X_t-\bar\xi\|_\infty \ge\vv\},\\
&\tau_n=\inf\{t\ge t_0: |X^i(t)-\bar X^i(t)|\ge n^{-1}\},\ \ n\ge 1.\end{split}\end{equation*}
Let $g_n(s)= \e^{ns}-1.$  Since $X^i_s\le \bar X^i_s$, and from (a) it follows that $$X^i(s)-\bar X^i(s) +\gg^i(s,X_s,z)- \bar\gg^i(s,\bar X_s,z)\le 0,\ \ s\ge t_0,$$ hence, by   It\^o's formula  we obtain
\beg{equation*}\beg{split} &0   \ge \E  g_n((X^i -\bar X^i)((t_0+\vv)\land \tt\tau\land \tau\land\tau_n))\\
& = \E \int_{t_0}^{(t_0+\vv)\land\tt\tau\land\tau\land\tau_n} \Big\{g_n'((X^i-\bar X^i)(s)) (b^i(s, X_s)-\bar b^i(s,\bar X_s))\\
&\ +\ff {g_n''(X^i(s)-\bar X^i(s))} 2\sum_{j=1}^m (\si^{ij}(s, X_s)-\bar\si^{ij}(s,\bar X_s))^2\\
&\  + \int_{E} \big\{g_n(X^i(s)-\bar X^i(s) +\gg^i(s,X_s,z)- \bar\gg^i(s,\bar X_s,z))-g_n(X^i(s)-\bar X^i(s))\big\}\nu(\d z)\Big\}\d s\\
&\ge \Big(\ff{n^2\vv^2}{2\e} -  Cn\e\Big)\E\{\vv\land(\tt\tau-t_0)\land(\tau-t_0)\land(\tau_n-t_0)\},\ \ n\ge 1,\end{split}\end{equation*} where, according to {\bf (A1)} and {\bf (A2)},
\beg{equation*}\beg{split} C:= \sup\big\{&|b^i(t,\eta)-\bar b^i(t,\bar\eta)|+ \int_{E} |\gg^i(t,\eta,z)-\bar\gg^i(t,\bar\eta,z)|\nu(\d z):\\
 &\ t\in [t_0, t_0+\vv], \|\eta-\xi\|_\infty+\|\bar\eta-\bar\xi\|_\infty\le\vv\big\}<\infty.\end{split}\end{equation*}    This implies $\E((\tt\tau-t_0)\land(\tau-t_0)\land(\tau_n-t_0))=0$ for large $n$, which is impossible   since $\P(A_\vv)>0$ and, due to the right-continuity of the solutions, $\tt\tau\land \tau_n\land \tau >t_0$ holds on the set $A_\vv.$
\end{proof}

\section{Existence and uniqueness of solutions}

 When $N=0$ and $b,\si$ are deterministic, the following result is included in \cite[Theorem 4.2]{SWY}. The appearance of $N$ makes the solution    discontinuous, so that the argument in the proof of \cite[Theorem 4.2]{SWY} leading to the existence of weak solutions by proving the tightness of the approximating solutions is no longer valid. Moreover, since the coefficients are now random, the Yamada-Watanabe principle used there is invalid either.  Due to   {\bf (A1)} and {\bf (A2)}, the proof of the uniqueness and non-explosion is standard. To prove the existence, we approximate the original equation by those with Lipschitz coefficients  and  construct a strong solution to the original equation by proving that the approximating solutions form a Cauchy sequence under the topology of locally uniform convergence.

\beg{thm}\label{T3.1} Let $b,\si,\gg$ satisfy {\bf (A1)} and {\bf (A2)} with $\bar b=0,\bar\si=0$ and $\bar\gg=0$. Then, for any $s\ge 0$ and
$\F_s$-measurable $\C$-valued random variable $\xi$, the  equation
 $$\d X(t)=b(t,X_t)\d t+\si(t,X_t)\d B(t)+\int_E \gg(t,X_{t-},z)N(\d t,\d z),\ \ t\ge s, X_s=\xi,$$ has a unique solution which satisfies
$$\E\Big(\sup_{t\in [s-r_0,T]} |X(t)|^2\Big|\F_s\Big)<\infty\ \ \text{a.s.\ for\ each\ }T>s.$$  \end{thm}

\beg{proof} Without loss of generality, we only prove for $s=0$ and simply write $\E_0=\E(\cdot|\F_0).$

 (a) $\E_0\sup_{t\le T} |X(t)|^2<\infty\ \P$-a.s. for any $T>0.$ Let $X(t)$ be a solution to the equation. Let
$$\tau_n=\inf\{t\ge 0: |X(t)| \ge \|\xi\|_\infty+n\},\ \ n\ge 1.$$ By   It\^o's formula,  we have
\beg{equation*}\beg{split} \d |X(s)|^2 = &2\<\si(s, X_s)\d B(s), X(s)\> + \big(\|\si(s,X_s)\|_{HS}^2+2\<b(s, X_s),X(s)\>\big)\d s \\
&+\int_E(|X_{s-}+\gg(s,X_{s-},z)|^2-|X_{s-}|^2)N(\d s,\d z).\end{split}\end{equation*} Then, by the same technique used in the proof of Theorem \ref{T1.1},   for any $T>0$ we may find  a constant $C(T)>0$ such that for any $n\ge 1$, the process $\phi_n(t):= \E_0\sup_{s\le t} |X(s\land\tau_n)|^2$  satisfies $\P$-a.s.
$$ \E_0\phi_n(t) \le C(T)  + C(T) \int_0^{t} \E_0 \phi_n(s) u(\E_0\phi_n(s))\d s,\ \ t\in [0,T].$$
Let $G(s)=\int_1^s\ff 1 {ru(r)}\d r, s>0.$ By   Bihari's inequality we have $\P$-a.s.
$$\E_0 \phi_n(t)\le G^{-1}\big(G(C(T)+\|\xi\|_\infty)+C(T) t\big)<\infty,\ \ t\in [0,T].$$ Letting $n\uparrow\infty$, we conclude that $\tau_n\uparrow\infty$ and
$\E_0\sup_{t\le T} |X(t)|^2<\infty$.

(b) The uniqueness of the solution. Let $X(t)$ and $\tt X(t)$ be two solutions to the equation with the same initial data $X_0$.  Again using the technique applied in the proof of Theorem \ref{T1.1},  for any $T>0$ we may find  a constant $C(T)>0$ such that the process $\phi(t):= \sup_{s\le t} |X(s)-\tt X(s)|^2$ satisfies $\P$-a.s.
$$\E_0 \phi(t) \le C(T)\int_0^t \E_0\phi(s) u(\E_0\phi(s))\d s,\ \ t\in [0,T].$$ Since $\int_0^1\ff 1 {su(s)}\d s=\infty$, by   Bihari's inequality we conclude that $\P$-a.s. $\E_0\phi(t)=0$ for $t\in [0,T].$ This implies that $\P$-a.s. $X(t)=\tt X(t)$ for all $t\ge 0$ since $T>0$ is arbitrary.

(c) Existence of the solution for bounded $b,\si$ and $\theta:=\int_E \big(|\gg(\cdot,\cdot,z)|^2+|\gg(\cdot,\cdot,z)|\big)\nu(\d z).$  If $u\equiv 1$, i.e. the coefficients are Lipschitz continuous in $\xi\in\C$ with respect to the uniform norm, then the existence and uniqueness of the solutions can be proved by a standard  argument (cf. \cite{ST}). To prove the existence of the solution,  we approximate the coefficients by using Lipschitz ones as follows. Let $\mu$ be the distribution of the $\C$-valued random variable $\mathbf B$ with $\mathbf   B(s):=   \tt B(r_0+1+s), s\in [-r_0,0],$ where $\tt B(s)$ is a $d$-dimensional Brownian motion with $\tt B(0)=0.$ For any $n\ge 1$, let
\beg{equation*}\beg{split} &b_n(t,\xi)=  \int_\C b(t, \xi+n^{-1}\eta)\mu(\d\eta), \ \si_n(t,\xi)= \int_\C \si(t,\xi+n^{-1}\eta)\mu(\d \eta),\\
 &\gg_n(t,\xi,z)=\int_\C \gg(t,\xi+n^{-1} \eta,z)\mu(\d z),\ \ t\ge 0, \xi\in \C, z\in E.\end{split}\end{equation*} Since $b,\si$ and $\theta$ are bounded, applying \cite[Corollary 1.3]{BWY} for $\si=\ff 1 n I_{d\times d}, m=0, Z=b=0$ and $T=1+r_0$, we conclude that, for any $n\ge 1$,
\beg{equation*}\beg{split} &|b_n(t,\xi)-b_n(t,\eta)|^2 +\|\si_n(t,\xi)-\si_n(t,\eta)\|_{HS}^2 + \int_E |\gg_n(t,\xi,z)-\gg_n(t,\eta,z)|^2 \nu(\d z)\\
&+\bigg(\int_E |\gg_n(t,\xi,z)-\gg_n(t,\eta,z)| \nu(\d z)\bigg)^2\le K_n(t)\|\xi-\eta\|_\infty^2\end{split}\end{equation*} holds for some positive $K_n\in C([0,\infty))$. Therefore, the equation
\beq\label{ED} \d X^{(n)}(t)= b_n(t, X^{(n)}_t ) \d t +\si_n (t,X^{(n)}_t)\d B(t) + \int_E \gg_n(t, X^{(n)}_t,z)N(\d t,\d z)\end{equation} starting from  $X^{(n)}_0=\xi$ has a unique solution. Moreover, by  Jensen's inequality,
we see that {\bf (A1)} and {\bf (A2)} hold for $b_n,\si_n$ and $\gg_n$ uniformly in $n\ge 1.$

Next,   by {\bf (A1)}, we may find a positive function $K\in C([0,\infty))$ such that
\beg{equation*}\beg{split} &|b_n(t,\xi)- b_l(t,\xi)|^2+ \|\si_n(t,\xi)-\si_l(t,\xi)\|_{HS}^2+ \int_E|\gg_n(t,\xi,z)-\gg_l(t,\xi,z)|^2\nu(\d z)\\
&+ \bigg(\int_E |\gg_n(t,\xi,z)-\gg_l(t,\xi,z)| \nu(\d z)\bigg)^2\le K(t) \vv_{n,l},\end{split}\end{equation*}
where, according to $\mu(\|\cdot\|_\infty^2)<\infty$, $su(s)\le c(1+s)$ for some constant $c>0$ and $su(s)\to 0$ as $s\to 0$ so that
$$\vv_{n,l}:= \int_\C \|(n^{-1} -l^{-1})\eta\|_\infty ^2u(\|(n^{-1} -l^{-1})\eta\|_\infty^2)\mu(\d \eta)\to 0\ \ \text{as}\ n,l\to\infty.$$ Combining this with {\bf (A1)} we obtain
\beq\label{PD}\beg{split} &|b_n(t,\xi)-b_l(t,\eta)|^2+\|\si_n(t,\xi)-\si_l(t,\eta)\|_{HS}^2+\int_E|\gg_n(t,\xi,z)-\gg_l(t,\eta,z)|^2\nu(\d z)\\
&+\bigg(\int_E|\gg_n(t,\xi,z)-\gg_l(t,\eta,z)|\nu(\d z)\bigg)^2\\
&\le K(t)\vv_{n,l} +K(t) \|\xi-\eta\|_\infty^2 u(\|\xi-\eta\|_\infty^2),\ \ \ t\ge 0, \xi,\eta\in\C,\end{split}\end{equation} for some positive $K\in C([0,\infty)).$  Moreover, since
{\bf (A1)} and {\bf (A2)} hold for $b_n,\si_n$ and $\gg_n$ uniformly in $n$, by (a) we have $\P$-a.s.
\beq\label{UNI} \sup_{n\ge 1} \E_0 \sup_{t\le T} |X^{(n)}(t)|^2<\infty,\ \ T>0.\end{equation}

 Now, as in (a) and (b), by   It\^o's formula,   Burkholder-Davis-Gundy's inequality, {\bf (A1)} and {\bf (A2)} holding for $b_n,\si_n$ and $\gg_n$ uniformly in $n$,    Jensen's inequality and (\ref{PD}), for any $T>0$ we may find  a constant $C(T)>0$ such that the process $\phi_{n,l}(t):= \sup_{s\le t} |X^{(n)}(s)-  X^{(l)}(s)|^2$ satisfies $\P$-a.s.
$$\E_0 \phi_{n,l}(t) \le C(T)\int_0^t \E_0\phi_{n,l}(s) u(\E_0\phi_{n,l}(s))\d s + C(T)\vv(n,l) ,\ \ t\in [0,T].$$ Since
 $\vv(n,l)\to 0$ as $n,l\to\infty$, by   Bihari's inequality and $\int_0^1\ff 1 {su(s)}\d s=\infty$, we
obtain $\P$-a.s.
$$\lim_{n,l\to\infty} \E_0 \sup_{t\le T} |X^{(n)}(t)-X^{(l)}(t)|^2=0,\ \ T>0.$$ Therefore,
as $n\to\infty$, the process  $X^{(n)}$ converges locally uniformly to a process $X$, which  solves the first equation in (\ref{1.0}) according  to   {\bf (A1)}, (\ref{UNI}) and the facts that $su(s)\le c(1+s)$ for some constant $c>0$ and $su(s)\to 0$ as $s\to 0.$

(d) Existence of the solution for unbounded $b,\si$ and $\theta$. For any $n\ge 1$, let $\vec{n}=(n,\cdots,n)\in\R^d.$ Define
$$\aa_n(\xi)= (\xi\land\vec{n})\lor (-\vec{n}),\ \ n\ge 1, \xi\in\C.$$ Let
$$b_n(t,\xi)= b(t\land n, \aa_n(\xi)),\ \ \si_n(t,\xi)= \si(t\land n, \aa_n(\xi)),\ \ \gg_n(t,\xi,z)= \gg(t\land n, \aa_n(\xi),z).$$ Then  $b_n,\si_n$ and $\theta_n:=\int_E\big(|\gg_n(\cdot,\cdot,z)|^2+|\gg_n(\cdot,\cdot,z)|\big)\nu(\d z)$ are bounded. Thus, according to (a)-(c), the equation (\ref{ED}) with $X_0^{(n)}=\xi$ has a unique solution $X^{(n)}(t), t\ge 0$. Since for any $l\ge n\ge 1$ and $\xi\in \C$ with $\|\xi\|_\infty\le n$ we have
$$b_n(t,\xi)=b_l(t,\xi), \ \ \si_n(t,\xi)=\si_l(t,\xi),\ \ \gg_n(t,\xi,z)=\gg_l(t,\xi,z),\ \ t\in [0,n],$$ by   uniqueness,
one has $X^{(n)}(t)= X^{(l)}(t)$ for $t\le\tau_n$ where
$$\tau_n:=n\land \inf\{t\ge 0: \|X^{(n)}_t\|_\infty\ge n\}.$$ Moreover, as in (a), we can prove that $\tau_n\uparrow \infty$ as $n\uparrow \infty$. Therefore, $X(t):= X^{(n)}(t)$ if $t<\tau_n$ gives rise to a solution of the original equation.
\end{proof}

\paragraph{Acknowledgement.} We would like to thank the referee for very   helpful comments and corrections, and Dr. Chenggui Yuan for  valuable discussions.

\end{document}